\newfont{\Bbb}{msbm10 scaled\magstephalf}
\documentclass[reqno]{amsart}
\usepackage{amssymb}
\usepackage{amsfonts}
\usepackage{amsmath, amsthm, amscd, amsfonts}
\usepackage{cases}
\usepackage{collref}
\usepackage{cite}

\newtheorem{Lemma}{Lemma}
\newtheorem{Theorem}{Theorem}

\begin{document}
\title[Quasi-similarity of multiplication operator]{On quasi-similarity of multiplication operator on the weighted Bergman space in the unit ball}

\author[C. Chen, Y. Wang and Y.X. Liang] {Cui Chen$^*$, Ya Wang and Yu-Xia Liang}

\address{\newline Cui Chen, Department of Mathematics, Tianjin University of Finance and Economics, Tianjin 300222, P.R. China}
\email{chencui\_cc@126.com}


\address{\newline Ya Wang, Department of Mathematics, Tianjin University of Finance and Economics, Tianjin 300222, P.R. China}
\email{wangyasjxsy0802@163.com}

\address{\newline  Yu-Xia Liang, School of Mathematical Sciences, Tianjin Normal University, Tianjin 300378, P.R. China.}
\email{liangyx1986@126.com}



\date{}
\thanks{\noindent $^*$Corresponding author.}

\begin{abstract}
For $\alpha>-1$, let $A_\alpha^2(\mathbb{B}_N)$ be the weighted Bergman space on the unit ball $\mathbb{B}_N$ in $\mathbb{C}^N$. In this paper, we prove that the multiplication operator $M_{z^n}$ is quasi-similar to $\oplus_1^{\prod_{i=1}^N n_i}M_z$ on $A_\alpha^2(\mathbb{B}_N)$ for the multi-index $n=(n_1,n_2,\cdots,n_N)$.
\end{abstract}

\maketitle

\section{Introduction}

Similarity of operators is a weaker concept than unitary equivalence. It is well known that invariant subspace of an operator can be identified in terms of invariant subspaces of a similar operator. Over the years, there are some results characterizing the lattice of closed invariant subspaces by constructing similar operators, we refer the interested readers to the recent papers such as \cite{AM,AM1,CP}.

In previous years, the similarity between $M_{z^n}$ and $\oplus_1^n M_z$ acting on a Hilbert space is an active topic which has be concerned in lots of papers. In 2007, Jiang and Li (see \cite{JL}) first obtained that analytic Toeplitz operator $M_{B(z)}$ is similar to $\oplus_1^n M_z$ on the Bergman space if and only if $B(z)$ is an $n$-Blaschke product. Next, Li (see \cite{L1}) in 2009 proved that multiplication operator $M_{z^n}$ is similar to $\oplus_1^n M_z$ on the weighted Bergman space. And then, Jiang and Zheng in \cite{JZ} extended the main result in \cite{JL} to the weighted Bergman space. In 2011, Douglas and Kim in \cite{DK} investigated the reducing subspaces for an analytic multiplication operator $M_{z^n}$ on the Bergman space $A_\alpha^2(A_r)$ of the annulus $A_r$. The similarity of $M_{z^n}$ and $\oplus_1^n M_z$ also holds on the weighted Hardy space and Sobolev disk algebra, they are shown in \cite{AH} and \cite{LLL1}. For further results related, see \cite{JS,ZRF}. Moreover, in 2017, the unitary equivalence of analytic multipliers on Sobolev disk algebra was discussed in \cite{CQW}. But such a characterization is not always hold for all Hilbert spaces. In Recent paper, Li, Lan and Liu (see \cite{LLL}) proved that multiplication operator $M_{z^n}$ is not similar to $\oplus_1^n M_z$ on the Fock space, they are quasi-similar, actually.

However, all these results are considered in the one-dimensional space, such a characterization is much more complicated for the high-dimension case. Based on the works above, in this paper we are interested in the corresponding result on the weighted Bergman space in the unit ball.

Let $\mathbb{B}_N$ be the unit ball in $\mathbb{C}^N$, where $\mathbb{C}^N$ is the $N$-dimensional complex vector space. Denoted by $\mathbb{N}_0$ and $\mathbb{N}$ the set of all positive and nonnegative integers, respectively. $H(\mathbb{B}_N)$ is the class of all holomorphic functions on $\mathbb{B}_N$.

Given three multi-indexes in $\mathbb{N}^N$ with $n=(n_1,n_2,\cdots,n_N), k=(k_1,k_2,\cdots,k_N)$ and $j=(j_1,j_2,\cdots,j_N)$, we use the notations
$$|n|=n_1+n_2+\cdots+n_N,\;\;\;\;\;\;n!=n_1!n_2!\cdots n_N!.$$
and
$$nk+j=(n_1k_1+j_1,n_2k_2+j_2,\cdots,n_Nk_N+j_N).$$
Moreover, we say $j<(\mbox{or}\leq)n$ if $j_i<(\mbox{or}\leq)n_i$ for each $i=1,2,\cdots,N$.

For $z\in\mathbb{C}^N$,
$$z^n=z_1^{n_1}z_2^{n_2}\cdots z_N^{n_N}.$$

Let $\alpha>-1$, the weighted Bergman space $A_\alpha^2(\mathbb{B}_N)$ consists of holomorphic functions $f$ in $L^2(\mathbb{B}_N,d v_\alpha)$, that is, $A_\alpha^2(\mathbb{B}_N)=L^2(\mathbb{B}_N,d v_\alpha)\cap H(\mathbb{B}_N)$. It is well known that $A_\alpha^2(\mathbb{B}_N)$ is a Hilbert space with the inner product defined by
$$\langle f,g\rangle=\int_{\mathbb{B}_N}f\overline{g}d v_\alpha,\;\;\;\;\;\;\mbox{for}\;f,g\in A_\alpha^2(\mathbb{B}_N).$$
The corresponding norm of $f$ is given by $\|f\|^2=\int_{\mathbb{B}_N}|f|^2d v_\alpha$.
As we know, the set of polynomials is dense in the Bergman space $A_\alpha^2(\mathbb{B}_N)$, and moreover,
$$\langle z^k,z^m\rangle=\int_{\mathbb{B}_N}z^k\overline{z^m}d v_\alpha=
\left\{\begin{array}{ll}
\frac{k!\Gamma(N+\alpha+1)}{\Gamma(N+|k|+\alpha+1)}\;,\;\;\;m=k,\\
\;\;\;\;\;\;\;\;0\;,\;\;\;\;\;\;\;\;\;\;\;\;\;\;otherwise.
\end{array}\right.$$
Thus $\{e_k(z)\triangleq\sqrt{\frac{\Gamma(N+|k|+\alpha+1)}{k!\Gamma(N+\alpha+1)}}z^k:k\in\mathbb{N}^N\}$ ia an orthonormal basis of $A_\alpha^2(\mathbb{B}_N)$.



Throughout the rest of this paper, we will always assume that $N\geq 2$ since it has been considered by Li \cite{L1} for the case $N=1$. Moreover, we would like to point out that $n$ stands for the multi-index $(n_1,n_2,\cdots,n_N)\in\mathbb{N}_0^N$ with $n_i\geq 2$ for at least one $i\in\{1,2,\cdots,N\}$, as it is trivial for otherwise.

\section{Quasi-similarity}

Let us recall that for two Hilbert spaces $\mathcal{H}$ and $\mathcal{K}$, an operator $X\in\mathcal{L}(\mathcal{H},\mathcal{K})$ is said to be quasi-invertible if it has zero kernel and dense range. Let $S\in\mathcal{L}(\mathcal{H})$ and $T\in\mathcal{L}(\mathcal{K})$, $S$ is quasi-similar to $T$ if there exist two quasi-invertible operators $X\in\mathcal{L}(\mathcal{H},\mathcal{K})$ and $Y\in\mathcal{L}(\mathcal{K},\mathcal{H})$ respectively such that $XS=TX$ and $SY=YT$. The following theorem is our main result.
\begin{Theorem}
$M_{z^n}$ is quasi-similar to $\oplus_1^{\prod_{i=1}^N n_i}M_z$ acting on the weighted Bergman space $A_\alpha^2(\mathbb{B}_N)$ with $\alpha>-1$.
\end{Theorem}

Before going to the main theorem, we need the following lemma.
\begin{Lemma}
For a multi-index $j=(j_1,j_2,\cdots,j_N), j_i=0,1,\cdots,n_i-1 \;\mbox{for each}\;i=1,2,\cdots,N$, let $\mathcal{A}_j=span\{e_{nk+j}:k=(k_1,k_2,\cdots,k_N)\in\mathbb{N}^N\}$. Then

(i) $\{e_{nk+j}:k=(k_1,k_2,\cdots,k_N)\in\mathbb{N}^N\}$ form an orthonormal basis of $\mathcal{A}_j$.

(ii) $A_\alpha^2(\mathbb{B}_N)=\oplus_{j}\mathcal{A}_j$, where the direct sum is over all multi-indexes $j=(j_1,j_2,\cdots,j_N)$ and $j_i=0,1,\cdots,n_i-1$ for each $i=1,2,\cdots,N$.

(iii) $\mathcal{A}_j$ is a reducing subspace for $M_{z^n}$.
\end{Lemma}
\textbf{Proof.}
(i) Note that for all $k, m\in\mathbb{N}^N$,
\begin{eqnarray*}
&&\langle e_{nk+j},e_{nm+j}\rangle\nonumber\\
&=&\int_{\mathbb{B}}\sqrt{\frac{\Gamma(N+|nk+j|+\alpha+1)}{(nk+j)!\Gamma(N+\alpha+1)}}z^{nk+j}
\sqrt{\frac{\Gamma(N+|nm+j|+\alpha+1)}{(nm+j)!\Gamma(N+\alpha+1)}}\overline{z}^{nm+j}dv_\alpha(z).
\end{eqnarray*}

If $m\neq k$, then $\langle e_{nk+j},e_{nm+j}\rangle=0$. If $m=k$, we have
\begin{eqnarray*}
\langle e_{nk+j},e_{nm+j}\rangle&=&\frac{\Gamma(N+|nk+j|+\alpha+1)}{(nk+j)!\Gamma(N+\alpha+1)}\int_{\mathbb{B}_N}|z^{nk+j}|^2dv_\alpha(z)\\
&=&\frac{\Gamma(N+|nk+j|+\alpha+1)}{(nk+j)!\Gamma(N+\alpha+1)}\frac{(nk+j)!\Gamma(N+\alpha+1)}{\Gamma(N+|nk+j|+\alpha+1)}=1.
\end{eqnarray*}
Thus (i) holds.

(ii) It is clear that $\mathcal{A}_j\perp\mathcal{A}_t$ for all $0\leq j\neq t\leq n-1$. Next, for $f\in A_\alpha^2(\mathbb{B}_N)$, it is easy to see that $f$ has the form
$$f=\sum_{k=0}^\infty\sum_{j=0}^{n-1}a_{jk}e_{nk+j}.$$
Suppose $f=0$, we conclude that $a_{jk}=0$ for all $k\in\mathbb{N}^N$ and $j=(j_1,j_2,\cdots,j_{N}),j_i=0,1,\cdots,n_i-1$, since
\begin{eqnarray*}
\langle\sum_{k=0}^\infty\sum_{j=0}^{n-1}a_{jk}e_{nk+j},e_l\rangle=0\;\;\;\;\;\;\mbox{for each}\;l\in\mathbb{N}^N.
\end{eqnarray*}
That is $0=\overbrace{0\oplus0\oplus\cdots\oplus0}^{\Pi_{i=1}^N n_i}$, which yields that $A_\alpha^2(\mathbb{B}_N)=\oplus_{j}\mathcal{A}_j$.

(iii) It is easy to see that both $\mathcal{A}_j$ and $\mathcal{A}_j^\perp$ are invariant subspaces for $M_{z^n}$.    $\Box$

By the previous lemma, it is clear that $M_{z^n}=\oplus_{j}M_{z^n}|_{\mathcal{A}_j}$. Then we can get the proof of Theorem 1.

\textbf{Proof of Theorem 1.}
Note that
$$M_ze_k=z\sqrt{\frac{\Gamma(N+|k|+\alpha+1)}{k!\Gamma(N+\alpha+1)}}z^k=\sqrt{\frac{(k+1)\Gamma(N+|k|+\alpha+1)}{\Gamma(N+|k+1|+\alpha+1)}}e_{k+1}.$$
Set $M_j=M_{z^n}|_{\mathcal{A}_j}$, then
\begin{eqnarray*}
M_je_{nk+j}&=&z^n\sqrt{\frac{\Gamma(N+|nk+j|+\alpha+1)}{(nk+j)!\Gamma(N+\alpha+1)}}z^{nk+j}\\
&=&\sqrt{\frac{(n(k+1)+j)!\Gamma(N+|nk+j|+\alpha+1)}{(nk+j)!\Gamma(N+|n(k+1)+j|+\alpha+1)}}e_{n(k+1)+j}
\end{eqnarray*}
Define $X_j:A_\alpha^2(\mathbb{B}_N)\rightarrow\mathcal{A}_j$ such that $X_je_k=c_{kj}e_{nk+j}$, where $c_{kj}$ are given by
\begin{eqnarray*}
c_{kj}=\sqrt{\frac{(nk+j)!\Gamma(N+|k|+\alpha+1)\Gamma(N+|j|+\alpha+1)}{k!j!\Gamma(N+|nk+j|+\alpha+1)\Gamma(N+\alpha+1)}}.
\end{eqnarray*}
Thus we conclude that $X_jM_ze_k=M_jX_je_k$. In fact,
\begin{eqnarray*}
X_j M_z e_k&=&X_j\sqrt{\frac{(k+1)\Gamma(N+|k|+\alpha+1)}{\Gamma(N+|k+1|+\alpha+1)}}e_{k+1}\\
&=&\sqrt{\frac{(nk+n+j)!\Gamma(N+|k|+\alpha+1)\Gamma(N+|j|+\alpha+1)}{k!j!\Gamma(N+|nk+n+j|+\alpha+1)\Gamma(N+\alpha+1)}}e_{nk+n+j},
\end{eqnarray*}
and
\begin{eqnarray*}
M_j X_j e_k&=&M_j\sqrt{\frac{(nk+j)!\Gamma(N+|k|+\alpha+1)\Gamma(N+|j|+\alpha+1)}{k!j!\Gamma(N+|nk+j|+\alpha+1)\Gamma(N+\alpha+1)}}e_{nk+j}\\
&=&\sqrt{\frac{(nk+n+j)!\Gamma(N+|k|+\alpha+1)\Gamma(N+|j|+\alpha+1)}{k!j!\Gamma(N+|nk+n+j|+\alpha+1)\Gamma(N+\alpha+1)}}e_{nk+n+j}.
\end{eqnarray*}
Now we will show that $X_j$ is not invertible in the following by giving the fact that $\liminf_{|k|\rightarrow\infty}c_{kj}=0$. Indeed,
\begin{eqnarray}
c_{kj}^2&=&\frac{\Gamma(N+|j|+\alpha+1)}{j!\Gamma(N+\alpha+1)}\frac{|nk+j|!\Gamma(N+|k|+\alpha+1)}{|k|!
\Gamma(N+|nk+j|+\alpha+1)}\frac{(nk+j)!|k|!}{|nk+j|!k!}\nonumber\\
&\triangleq&\frac{\Gamma(N+|j|+\alpha+1)}{j!\Gamma(N+\alpha+1)}I_1 I_2,\label{1}
\end{eqnarray}
where
$$I_1=\frac{|nk+j|!\Gamma(N+|k|+\alpha+1)}{|k|!\Gamma(N+|nk+j|+\alpha+1)}\;\;\;\;\;\mbox{and}\;\;\;\;\;I_2=\frac{(nk+j)!|k|!}{|nk+j|!k!}.$$
For $I_1$,
\begin{eqnarray}
&&\frac{|nk+j|!\Gamma(N+|k|+\alpha+1)}{|k|!\Gamma(N+|nk+j|+\alpha+1)}\nonumber\\
&=&\frac{|nk+j|(|nk+j|-1)\cdots(|k|+1)}{(N+|nk+j|+\alpha)(N+|nk+j|+\alpha-1)\cdots(N+|k|+\alpha+1)}\nonumber\\
&=&\frac{|k|+1}{N+|nk+j|+\alpha}\frac{(|nk+j|)(|nk+j|-1)\cdots(|k|+2)}{(N+|nk+j|+\alpha-1))\cdots(N+|k|+\alpha+1)}\nonumber\\
&=&\frac{|k|+1}{N+|nk+j|+\alpha}\frac{1}{(1+\frac{N+\alpha-1}{|nk+j|})(1+\frac{N+\alpha-1}{|nk+j|-1})\cdots(1+\frac{N+\alpha-1}{|k|+2})}\label{2}
\end{eqnarray}
Define
\begin{eqnarray*}
a_{kj}(\alpha)=(1+\frac{N+\alpha-1}{|nk+j|})(1+\frac{N+\alpha-1}{|nk+j|-1})\cdots(1+\frac{N+\alpha-1}{|k|+2}).
\end{eqnarray*}
Since $\alpha>-1$ and $N\geq 2$, it is clear $N+\alpha-1>0$. Then it easily follows that
\begin{eqnarray}
a_{kj}(\alpha)&\leq&(1+\frac{N+\alpha-1}{|k|+2})^{|nk+j|-|k|-1}\nonumber\\
&\leq&(1+\frac{N+\alpha-1}{|k|+2})^{|n||k|+|j|-|k|-1}\nonumber\\
&=&(1+\frac{N+\alpha-1}{|k|+2})^{\frac{|k|+2}{N+\alpha-1}\cdot\frac{N+\alpha-1}{|k|+2}\cdot[|n||k|+|j|-|k|-1]}\nonumber\\
&\rightarrow&\exp\{(|n|-1)(N+\alpha-1)\},\;\;\;\mbox{as}\;|k|\rightarrow\infty.\label{3}
\end{eqnarray}
On the other hand,
\begin{eqnarray}
a_{kj}(\alpha)&\geq&(1+\frac{N+\alpha-1}{|nk+j|})^{|nk+j|-|k|-1}\nonumber\\
&\geq&(1+\frac{N+\alpha-1}{|n||k|+|j|})^{|k|+|j|-|k|-1}\nonumber\\
&=&(1+\frac{N+\alpha-1}{|n||k|+|j|})^{|j|-1}\nonumber\\
&\rightarrow&1,\;\;\;\mbox{as}\;|k|\rightarrow\infty.\label{4}
\end{eqnarray}
From (\ref{2})-(\ref{4}), associated with
\begin{eqnarray*}
\frac{1}{|n|}&=&\lim_{|k|\rightarrow\infty}\frac{|k|+1}{N+|n||k|+|j|+\alpha}\\
&\leq&\lim_{|k|\rightarrow\infty}\frac{|k|+1}{N+|nk+j|+\alpha}\\
&\leq&\lim_{|k|\rightarrow\infty}\frac{|k|+1}{N+n_m|k|+\alpha}=\frac{1}{n_m},
\end{eqnarray*}
where $m\in\{1,2,\cdots,N\}$ is the index such that $n_m=\min\{n_1,n_2,\cdots,n_N\}$, then we conclude that there exist two positive constants $C_1$ and $C_2$ such that
\begin{eqnarray}
C_1\leq\liminf_{|k|\rightarrow\infty}I_1\leq\limsup_{|k|\rightarrow\infty}I_1\leq C_2.\label{5}
\end{eqnarray}
For $I_2$, first note that
\begin{eqnarray}
\frac{(n_1k_1+j_1)!\cdots(n_N k_N+j_N)!}{(n_1k_1+j_1+\cdots+n_N k_N+j_N)!}
\leq\frac{(n_m k_1+j_m)!\cdots(n_m k_N+j_m)!}{(n_m k_1+j_m+\cdots+n_m k_N+j_m)!}.\label{6}
\end{eqnarray}
This is due to that
\begin{eqnarray*}
&&\frac{(n_1k_1+j_1)!\cdots(n_N k_N+j_N)!}{(n_m k_1+j_m)!\cdots(n_m k_N+j_m)!}
\cdot\frac{(n_m k_1+j_m+\cdots+n_m k_N+j_m)!}{(n_1k_1+j_1+\cdots+n_N k_N+j_N)!}\\
&=&\frac{\Pi_{i=1}^N[(n_i k_i+j_i)(n_i k_i+j_i-1)\cdots(n_m k_i+j_m+1)]}
{\sum_{i=1}^N(n_i k_i+j_i)\{\sum_{i=1}^N(n_i k_i+j_i)-1\}\cdots\{\sum_{i=1}^N(n_m k_i+j_m)+1\}}\leq 1,
\end{eqnarray*}
which comes from the numbers of factors are both $\sum_{i=1}^N(n_i k_i+j_i-n_m k_i-j_m)$ in the numerator and denominator, and the molecular factors are not greater than the denominator factors one by one.

Using Stirling's approximation, we have
\begin{eqnarray*}
&&\liminf_{|k|\rightarrow\infty}I_2\leq\lim_{k_1=k_2=\cdots=k_N\rightarrow\infty}I_2\\
&=&\lim_{k_1\rightarrow\infty}\frac{[(n_m k_1+j_m)!]^N}{[N(n_m k_1+j_m)]!}\frac{(N k_1)!}{(k_1!)^N}\\
&=&\lim_{k_1\rightarrow\infty}\frac{[(n_m k_1+j_m)^{n_m k_1+j_m}\sqrt{n_m k_1+j_m}e^{-(n_m k_1+j_m)}]^N}
{[N(n_m k_1+j_m)]^{N(n_m k_1+j_m)}\sqrt{N(n_m k_1+j_m)}e^{-N(n_m k_1+j_m)}}\\
&&\cdot\frac{(N k_1)^{N k_1}\sqrt{N k_1}e^{-N k_1}}{[k_1^{k_1}\sqrt{k_1}e^{-k_1}]^N}\\
&=&\lim_{k_1\rightarrow\infty}n_m^{(N-1)/2}(N^N)^{(1-n_m)k_1-j_m}
\end{eqnarray*}
Thus if $n_m>1$, $\liminf_{|k|\rightarrow\infty}I_2=0$.

For the remaining case $n_m=1$. Without loss of generality, we may assume that $\sharp\{i\in\{1,2,\cdots,N\}: n_i=1\}=1$, for otherwise set $k_{m'}=0$ when $n_{m'}=1$. Choose $s$ such that $n_s=\min\{n_i:i\in\{1,2,\cdots,N\}\backslash\{m\}\}$. From the same reason of (\ref{6}) we can get
\begin{eqnarray*}
&&\frac{(n_1k_1+j_1)!\cdots(n_N k_N+j_N)!}{(n_1k_1+j_1+\cdots+n_N k_N+j_N)!}\\
&\leq&\frac{(n_s k_1+j_s)!\cdots(n_s k_{m-1}+j_s)!(k_m+j_m)!(n_s k_{m+1}+j_s)!\cdots(n_s k_N+j_s)!}{(n_s k_1+j_s+\cdots+n_s k_{m-1}+j_s+k_m+j_m+n_s k_{m+1}+j_s+\cdots+n_s k_N+j_s)!}.
\end{eqnarray*}
By Stirling's approximation, we conclude that
\begin{eqnarray*}
&&\liminf_{|k|\rightarrow\infty}I_2\leq\lim_{k_1=k_2=\cdots=k_N\rightarrow\infty}I_2\\
&=&\lim_{k_1\rightarrow\infty}\frac{[(n_s k_1+j_s)!]^{N-1}(k_1+j_m)!}{[(N-1)(n_s k_1+j_s)+k_1+j_m]!}\frac{(N k_1)!}{(k_1!)^N}\\
&=&\lim_{k_1\rightarrow\infty}\frac{(n_s k_1+j_s)^{(N-1)(n_s k_1+j_s)}\sqrt{(n_s k_1+j_s)^{N-1}}e^{-(N-1)(n_s k_1+j_s)}}
{[(N-1)(n_s k_1+j_s)+k_1+j_m]^{(N-1)(n_s k_1+j_s)+k_1+j_m}}\\
&&\cdot\frac{(k_1+j_m)^{k_1+j_m}\sqrt{k_1+j_m}e^{-(k_1+j_m)}}{\sqrt{(N-1)(n_s k_1+j_s)+k_1+j_m}e^{-[(N-1)(n_s k_1+j_s)+k_1+j_m]}}\\
&&\cdot\frac{(N k_1)^{N k_1}\sqrt{N k_1}e^{-N k_1}}{(k_1^{k_1}\sqrt{k_1}e^{-k_1})^N}\\
&=&\lim_{k_1\rightarrow\infty}\frac{N^{N k_1}}
{(N-1+\frac{k_1+j_m}{n_s k_1+j_s})^{(N-1)(n_s k_1+j_s)}(1+\frac{(N-1)(n_s k_1+j_s)}{k_1+j_m})^{k_1+j_m}}\\
&&\cdot\frac{\sqrt{(n_s k_1+j_s)^{N-1}}\sqrt{k_1+j_m}\sqrt{N k_1}}{\sqrt{(N-1)(n_s k_1+j_s)+k_1+j_m}\sqrt{k_1^N}}\\
&=&\lim_{k_1\rightarrow\infty}\Big(\frac{N}{(N-1+\frac{1}{n_s})^{n_s}}\Big)^{(N-1)k_1}
\Big(\frac{N}{1+(N-1)n_s}\Big)^{k_1}\sqrt{\frac{N n_s^{N-1}}{(N-1)n_s+1}}\\
&&\cdot\frac{1}{(N-1+\frac{1}{n_s})^{j_s}(1+(N-1)n_s)^{j_m}}\\
&=&0,
\end{eqnarray*}
which we have used $\Big(N-1+\frac{1}{n_s}\Big)^{n_s}>N$ and $1+(N-1)n_s>N$ in the last equation, since $N, n_s\geq 2$. Along with (\ref{1}) and (\ref{5}), we have $\liminf_{|k|\rightarrow\infty}c_{k,j}=0$. This implies that $X_j$ is not invertible.

Next, define $Y_j:\mathcal{A}_j\rightarrow A_\alpha^2(\mathbb{B}_N)$ by $Y_je_{nk+j}=b_{kj}e_k$, where $b_{kj}=\frac{1}{c_{kj}}$. It is easy to prove that $Y_j M_j=M_z Y_j$ on $\mathcal{A}_j$. From the above discussion, we have $\limsup_{|k|\rightarrow\infty}b_{kj}=\infty$. Hence, $Y_j$ is also not invertible.

For the rest of the proof we will show that both $X_j, Y_j$ are quasi-invertible operators. To do this, take $f\in\mbox{ker}X_j$ with $f=\sum_{k=0}^\infty d_k e_k, d_k\in\mathbb{C}^N$. From $0=\langle X_jf,e_{nk+j}\rangle=\langle\sum_{k=0}^\infty d_k c_{kj} e_{nk+j},e_{nl+j}\rangle$, we have $d_l=0$ for all $l\in\mathbb{N}^N$, Hence $X_j$ is injective. Besides, pick $g\in\mbox{ker}X_j^*$ with $g=\sum_{k=0}^\infty m_k e_{nk+j}, m_k\in\mathbb{C}^N$. From $0=\langle e_l,X_j^*g\rangle=\langle c_{lj}e_{nl+j},\sum_{k=0}^\infty m_k e_{nk+j}\rangle$, we have $m_l=0$ for all $l\in\mathbb{N}^N$. Then $(\mbox{Ran}X_j)^\perp=\mbox{ker}X_j^*=\{0\}$, so $\overline{\mbox{Ran}X_j}=\mathcal{A}_j$, that is $X_j$ has dense range. Similarly, we can check that $Y_j$ is injective and has dense range. It turns out that both $X_j$ and $Y_j$ are quasi-invertible. To sum up, $M_j$ is quasi-similar to $M_z$. Since $M_{z^n}=\oplus_j M_j$, we conclude that $M_{z^n}$ is quasi-similar to $\oplus_1^{\prod_{i=1}^N n_i}M_z$, and we are done. $\;\;\;\;\;\Box$

\end{document}